\documentclass{article}
\usepackage{xcolor}
\usepackage{amsmath}
\usepackage{amssymb}
\usepackage{amsfonts}
\usepackage{amsthm}
\usepackage{blkarray}
\usepackage{arydshln}
\usepackage{multirow}
\newcommand{\bigzerou}{%
\smash{\lower1.7ex\hbox{\bg 0}}}
\theoremstyle{definition}
\newtheorem{theorem}{Theorem}[section]

\newtheorem{defi}{Definition}[section]

\newtheorem{conj}{Conjecture}[section]

\newtheorem{ex}{Example}[section]

\newtheorem{err}{Erratum}[section]

\newcommand{\ba}{\begin{eqnarray}}
\newcommand{\ea}{\end{eqnarray}}
\newcommand{\ban}{\begin{eqnarray*}}
\newcommand{\ean}{\end{eqnarray*}}
\newcommand{\no}{\nonumber}

\newcommand{\mapright}[1]{%
\smash{\mathop{%
\hbox to 1.0cm{\rightarrowfill}}\limits^{#1}}}
\newcommand{\mapleft}[1]{%
\smash{\mathop{%
\hbox to 1.3cm{\leftarrowfill}}\limits^{#1}}}
\begin{document}
\title{\begin{flushright}
\begin{minipage}[b]{5em}
\normalsize
${}$ \\
\end{minipage}
\end{flushright}
{\bf Geometrical Proof of Generalized Mirror Transformation for Multi-Point Virtual Structure Constants of Projective Hypersurfaces}}
\author{Masao Jinzenji\\
\\
\it Department of Mathematics, Okayama University\\
\it Okayama, 700-8530, Japan\\
{\it e-mail address: pcj70e4e@okayama-u.ac.jp}}
\maketitle

\begin{abstract}
In this paper, we propose a geometric proof of the generalized mirror transformation for multi-point virtual structure constants of degree $k$ hypersurfaces in $CP^{N-1}$.
\end{abstract}
\section{Introduction}
\subsection{Notation and Main Theorem}
Let $N$ and $k$ be positive integers, and let $M_{N}^{k}$ denote a degree $k$ hypersurface in ${CP}^{N-1}$. We consider $\overline{M}_{0,n}(CP^{N-1},d)$, the moduli space of stable maps of degree $d$ from genus 0 semi-stable curves with $n$ marked points to ${CP}^{N-1}$ \cite{K}. The genus 0 $n$-pointed Gromov-Witten invariant of $M_{N}^{k}$ is defined as \cite{K,Jin2}:
 \ba
\langle\prod_{j=1}^{n}{\cal O}_{h^{c_{j}}}\rangle_{0,d}=\int_{\overline{M}_{0,n}(CP^{N-1},d)}
c_{top}(\overline{\cal E}_{d}^{k})\wedge\biggl(\mathop{\bigwedge}_{j=1}^{n} ev_{j}^{*}(h^{c_{j}})\biggr),
\ea
where $h$ is the hyperplane class in $H^{*}({CP}^{N-1},{\bf C})$, and $ev_{j}: \overline{M}_{0,n}({CP}^{N-1},d) \rightarrow {CP}^{N-1}$ ($j=1,2,\cdots,n$) are the evaluation maps at the $j$-th marked point. Here, $\overline{\cal E}_{d}^{k}$ is the vector bundle on $\overline{M}_{0,n}({CP}^{N-1},d)$ that imposes the condition that the image of the stable map is contained in the hypersurface. This invariant is non-vanishing only if the following condition holds:
\ba
\sum_{j=1}^{n}c_{j}=N-5+(N-k)d+n.
\ea

We also introduce the compactified moduli space of quasimaps from ${CP}^{1}$ with $2+n$ marked points ($0$, $\infty$, and $z_{j} \in {CP}^{1} - \{0, \infty\}$ for $j=1,2,\cdots,n$) to ${CP}^{N-1}$ of degree $d$, denoted by $\widetilde{Mp}_{0,2|n}(N,d)$. We have to note here that the latter $n$ marked points do not have to be distinct. For details on the construction, see \cite{Jin1} for the case of two marked points, and \cite{JS} for the case of $2+n$ marked points.   We define the intersection number $w({\cal O}_{h^{a}}{\cal O}_{h^{b}}|\prod_{j=1}^{n}{\cal O}_{h^{c_{j}}})_{0,d}$ on $\widetilde{Mp}_{0,2|n}(N,d)$, analogous to $\langle{\cal O}_{h^{a}}{\cal O}_{h^{b}} \prod_{j=1}^{n}{\cal O}_{h^{c_{j}}}\rangle_{0,d}$ on $\overline{M}_{0,2+n}(\mathbb{CP}^{N-1},d)$, as:
\ba
w({\cal O}_{h^{a}}{\cal O}_{h^{b}}|\prod_{j=1}^{n}{\cal O}_{h^{c_{j}}})_{0,d}=\int_{\widetilde{Mp}_{0,2|n}(N,d)}
c_{top}(\widetilde{\cal E}_{d}^{k})\wedge ev_{0}^{*}(h^a)\wedge ev_{\infty}^{*}(h^b)\wedge\left(\mathop{\bigwedge}_{j=1}^{n}ev_{j}^{*}(h^{c_{j}})\right).
\label{defw}
\ea  
Here, $ev_{0}$ and $ev_{\infty}$ are the evaluation maps at $0$ and $\infty$, respectively, $ev_{j}$ is the evaluation map at $z_{j}$, and $\widetilde{\cal E}_{d}^{k}$ is the vector bundle on $\widetilde{Mp}_{0,2|n}(N,d)$ with the same geometric meaning as $\overline{\cal E}_{d}^{k}$. This intersection number is non-vanishing only if the following condition is satisfied:
\ba
a+b+\sum_{j=1}^{n}(c_{j}-1)=N-3+(N-k)d\;\;\left(\Longleftrightarrow\;\; a+b+\sum_{j=1}^{n}c_{j}=N-3+(N-k)d+n\right).
\label{wsel}
\ea
We refer to this intersection number as the ``multi-point virtual structure constant'' \cite{JS}.

Let $P_{g}^{l}$ be the set of partitions of the positive integer $g$ into $l$ parts:
\ba
P_{g}^{l}=\{\sigma_{g}=(g_{1},\cdots,g_{l})\;|\;1\leq g_{1}\leq g_{2}\leq \cdots\leq g_{l},\;\;
\sum_{j=1}^{l}g_{j}=g\;\}.
\ea
For a partition $\sigma_{g} \in P_{g}^{l}$, we define the multiplicity $\mbox{mul}(i,\sigma_{g})$ of $\sigma_{g}$ as:
\ba
\mbox{mul}(i,\sigma_{g})=(\mbox{number of subscripts $j$ such that $g_{j}=i$}).
\ea
We define the combinatorial factor $S(\sigma_{g})$ as:
\ba
S(\sigma_{g})=\prod_{i=1}^{g}\frac{1}{(\mbox{mul}(i,\sigma_{g}))!}.
\ea

In this paper, we prove the following theorem, which describes the relationship between the intersection numbers $w({\cal O}_{h^{a}}{\cal O}_{h^{b}}|\prod_{j=1}^{n}{\cal O}_{h^{c_{j}}})_{0,d}$ and the Gromov-Witten invariants $\langle\prod_{j=1}^{n}{\cal O}_{h^{c_{j}}}\rangle_{0,d}$.
\begin{theorem}
The following equality holds:
\begin{eqnarray}
&&w({\cal O}_{h^{a}}{\cal O}_{h^{b}}
|\prod_{j=2}^{N-2}({\cal O}_{h^{j}})^{n_{j}})_{0,d}\no\\
&&=\langle{\cal O}_{h^{a}}{\cal O}_{h^{b}}\prod_{j=2}^{N-2}({\cal O}_{h^{j}})^{n_{j}}\rangle_{0,d}+w({\cal O}_{h^{a+b}}{\cal O}_{1}|\prod_{j=2}^{N-2}({\cal O}_{h^{j}})^{n_{j}})_{0,d}+ \no\\
&&+\sum_{g=1}^{d-1}\sum_{l=1}^{g}\sum_{\substack{\sigma_{g}\in P_{g}^{l}\\(\sigma_{g}=(g_{1},\cdots,g_{l}))}}
S(\sigma_{g})\Biggl(\sum_{\substack{m_{0}^{j},m_{1}^{j},m_{2}^{j},\cdots,m_{l}^{j}\geq 0\\m_{0}^{j}+\sum_{i=1}^{l}m_{i}^{j}=n_{j}}}\sum_{e_{1},\cdots,e_{l}=0}^{N-2}
\left(\prod_{j=2}^{N-2}\frac{n_{j}!}{m_{0}^{j}!\prod_{p=1}^{l}m_{p}^{j}!}\right)\no\\
&&\times
\langle{\cal O}_{h^{a}}{\cal O}_{h^{b}}\prod_{j=2}^{N-2}({\cal O}_{h^{j}})^{m_{0}^{j}}\prod_{p=1}^{l}{\cal O}_{h^{e_{p}}}\rangle_{0,d-g}
\biggl(\prod_{p=1}^{l}\frac{w({\cal O}_{h^{N-2-e_{p}}}{\cal O}_{1}|\prod_{j=2}^{l}({\cal O}_{h^{j}})^{m_{p}^{j}})_{0,g_{p}}}{k}\biggr)\Biggr),
\label{main}
\end{eqnarray}
where
\ban
a+b+\sum_{j=2}^{N-2}(j-1)n_{j}=N-3+(N-k)d.
\ean
\label{gt}
\end{theorem}

\subsection{Usage and Examples of the Main Theorem}
In \cite{JS}, we established the following properties for multi-point virtual structure constants:
\begin{itemize}
\item[(i)]{$w({\cal O}_{h^{a}}{\cal O}_{h^{b}}|{\cal O}_{h^{c}})_{0,0}=k\cdot\delta_{N-2,a+b+c}$.}
\item[(ii)]{If $n\neq 1$, then $w({\cal O}_{h^{a}}{\cal O}_{h^{b}}|\prod_{j=1}^{n}{\cal O}_{h^{c_{j}}})_{0,0}=0$.}
\item[(iii)]{$w({\cal O}_{h^{a}}{\cal O}_{h^{b}}|\prod_{p=0}^{N-2}({\cal O}_{h^{p}})^{m_{p}})_{0,d}=\delta_{m_{0},0}\cdot d^{m_{1}}w({\cal O}_{h^{a}}{\cal O}_{h^{b}}|\prod_{p=2}^{N-2}({\cal O}_{h^{p}})^{m_{p}})_{0,d}$ for $d\geq 1$.}
\end{itemize}
The characteristics (i) and (ii) follow from the fact that the degrees of freedom of $n(\geq 2)$ marked points cannot be fixed by insertion of operators in the $d=0$ case. The property (iii) follows from 
Theorem 1 in \cite{JS}. 
Consequently, we will only consider $w({\cal O}_{h^{a}}{\cal O}_{h^{b}}|\prod_{p=2}^{N-2}({\cal O}_{h^{p}})^{m_{p}})_{0,d}$ for $d\geq 1$ henceforth.
Let us introduce the generating function for the multi-point virtual structure constants:
\ba
&&w({\cal O}_{h^a}{\cal O}_{h^b})_{0}(x^{0},x^{1},\cdots,x^{N-2})\no\\
&=&kx^{N-2-a-b}+\sum_{m_{2}=0}^{\infty}\cdots\sum_{m_{N-2}=0}^{\infty}\sum_{d=1}^{\infty}e^{dx^{1}}w({\cal O}_{h^{a}}{\cal O}_{h^{b}}|\prod_{p=2}^{N-2}({\cal O}_{h^{p}})^{m_{p}})_{0,d}\prod_{q=2}^{N-2}\frac{(x^{q})^{m_q}}{m_{q}!}.
\label{simp}
\ea
In \cite{JS}, we defined the ``mirror map'', a coordinate change of deformation variables, as:
\ba
t^p(x^{0},x^{1},\cdots,x^{N-2}):=\frac{1}{k}w({\cal O}_{h^{N-2-p}}{\cal O}_{1})_{0}\;\;\;(p=0,1,\cdots, N-2).
\label{mmap}
\ea 
Equation (\ref{simp}) reveals the structure of this mirror map:
\ba
t^p(x^{0},x^{1},\cdots,x^{N-2})=
x^{p}+\frac{1}{k}\sum_{m_{2}=0}^{\infty}\cdots\sum_{m_{N-2}=0}^{\infty}\sum_{d=1}^{\infty}e^{dx^{1}}w({\cal O}_{h^{N-2-p}}{\cal O}_{1}|\prod_{p=2}^{N-2}({\cal O}_{h^{p}})^{m_{p}})_{0,d}\prod_{q=2}^{N-2}\frac{(x^{q})^{m_q}}{m_{q}!}.
\label{mmap2}
\ea
This structure allows us to invert the mirror map, yielding:
\ba
x^{p}=x^{p}(t^{0},t^{1},\cdots,t^{N-2})\;\;\;(p=0,1,\cdots,N-2).
\label{invert}
\ea

On the other hand, Gromov-Witten invariants satisfy the following properties \cite{KM}:
\begin{itemize}
\item[(i)]{$\langle {\cal O}_{h^{a}}{\cal O}_{h^{b}}{\cal O}_{h^{c}}\rangle_{0,0}=k\cdot\delta_{N-2,a+b+c}.$}
\item[(ii)] {If $d=0$ and $n\neq 3$, $\langle\prod_{j=1}^{n}{\cal O}_{h^{c_{j}}}\rangle_{0,0}$ vanishes.}  
\item[(iii)] {$\langle \prod_{p=0}^{N-2}({\cal O}_{h^{p}})^{m_{p}}\rangle_{0,d}=\delta_{m_{0},0}\cdot d^{m_{1}}\langle \prod_{p=2}^{N-2}({\cal O}_{h^{p}})^{m_{p}}\rangle_{0,d}\quad(d\geq 1).
$}
\end{itemize}
Analogous to the generating function in (\ref{simp}), we introduce the perturbed two-point genus 0 Gromov-Witten invariant:
\ba
&&\langle{\cal O}_{h^a}{\cal O}_{h^b}\rangle_{0}(t^{0},t^{1},\cdots,t^{N-2})\no\\
&:=&\sum_{m_{0}=0}^{\infty}\cdots\sum_{m_{N-2}=0}^{\infty}\sum_{d=0}^{\infty}\langle{\cal O}_{h^{a}}{\cal O}_{h^{b}}\prod_{p=0}^{N-2}({\cal O}_{h^{p}})^{m_{p}}\rangle_{0,d}
\prod_{q=0}^{N-2}\frac{(t^{q})^{m_q}}{m_{q}!}\no\\
&=&kt^{N-2-a-b}+\sum_{m_{2}=0}^{\infty}\cdots\sum_{m_{N-2}=0}^{\infty}\sum_{d=1}^{\infty}e^{dt^{1}}\langle{\cal O}_{h^{a}}{\cal O}_{h^{b}}\prod_{p=2}^{N-2}({\cal O}_{h^{p}})^{m_{p}}\rangle_{0,d}\prod_{q=2}^{N-2}\frac{(t^{q})^{m_q}}{m_{q}!}.
\ea
Then, the conjecture proposed in \cite{JS} is:

\begin{conj}{\bf \cite{JS}}
\ba
&&\langle{\cal O}_{h^a}{\cal O}_{h^b}\rangle_{0}(t^{0},t^{1},\cdots,t^{N-2})=w({\cal O}_{h^a}{\cal O}_{h^b})_{0}(x^{0}(t^{*}),x^{1}(t^{*}),\cdots,x^{N-2}(t^{*}))\no\\
&&(\Longleftrightarrow \langle{\cal O}_{h^a}{\cal O}_{h^b}\rangle_{0}(t^{0}(x^{*}),t^{1}(x^{*}),\cdots,t^{N-2}(x^{*}))=w({\cal O}_{h^a}{\cal O}_{h^b})_{0}(x^{0},x^{1},\cdots,x^{N-2})).
\label{jsconj}
\ea
Here, $x^{p}(t^{*})$ (resp. $t^{p}(x^{*})$) abbreviates $x^{p}(t^{0},t^{1},\cdots,t^{N-2})$ (resp. $t^{p}(x^{0},x^{1},\cdots,x^{N-2})$) as given in (\ref{invert}) (resp. (\ref{mmap})).
\label{jsc}
\end{conj} 
By expanding the right-hand side of the second line of (\ref{jsconj}) after substituting (\ref{mmap2}), we observe that Theorem \ref{gt} is equivalent to this conjecture. Therefore, this paper provides a proof of Conjecture \ref{jsc}.
\begin{ex}
\label{cp2}
Consider the case $N=4, k=1$, corresponding to the ${CP}^{2}$ model. As shown in \cite{Jin6}, the generalized mirror transformation for the two-point virtual structure constant is trivial:
\ba
w({\cal O}_{h^{2}}{\cal O}_{h^{2}})_{0,1}=\langle{\cal O}_{h^{2}}{\cal O}_{h^{2}}\rangle_{0,1}.
\ea
However, for multi-point virtual structure constants, the transformation can be non-trivial even for $d=1$. For $w({\cal O}_{h}{\cal O}_{h^{2}}|{\cal O}_{h^2})_{0,1}$, Theorem \ref{gt} yields:
\ba
w({\cal O}_{h}{\cal O}_{h^{2}}|{\cal O}_{h^2})_{0,1}=\langle{\cal O}_{h}{\cal O}_{h^{2}}{\cal O}_{h^{2}}\rangle_{0,1}
(=1)..
\label{exam2}
\ea
While this specific case is trivial, the theorem also indicates:
\ba
(2=)w({\cal O}_{h}{\cal O}_{h}|{\cal O}_{h^2}{\cal O}_{h^{2}})_{0,1}&=&\langle{\cal O}_{h}{\cal O}_{h}{\cal O}_{h^2}{\cal O}_{h^{2}}\rangle_{0,1}+
\langle{\cal O}_{h}{\cal O}_{h}{\cal O}_{1}\rangle_{0,0}w({\cal O}_{h^{2}}{\cal O}_{1}|{\cal O}_{h^{2}}{\cal O}_{h^2})_{0,1}\no\\
&=&\langle{\cal O}_{h}{\cal O}_{h}{\cal O}_{h^2}{\cal O}_{h^{2}}\rangle_{0,1}+w({\cal O}_{h^{2}}{\cal O}_{1}|{\cal O}_{h^{2}}{\cal O}_{h^2})_{0,1}(=1+1),
\label{exam20}
\ea 
which is non-trivial. The reason for this difference will be discussed in the next section.
\end{ex}

\begin{ex}
\label{cy88}
Next, we consider the $N=k=8$ case, representing a Calabi-Yau hypersurface. In this scenario, the generalized mirror transformation for the two-point virtual structure constants is non-trivial. For $d=1$, Theorem \ref{gt} provides the following equality:
\ba
(83871744=)w({\cal O}_{h^{2}}{\cal O}_{h^{2}}|{\cal O}_{h^{2}})_{0,1}&=&\langle{\cal O}_{h^{2}}{\cal O}_{h^{2}}{\cal O}_{h^{2}}\rangle_{0,1}+\langle{\cal O}_{h^{2}}{\cal O}_{h^{2}}{\cal O}_{h^{2}}\rangle_{0,0}w({\cal O}_{h^{4}}{\cal O}_{1}|{\cal O}_{h^{2}})_{0,1}/k\no\\
&=&\langle{\cal O}_{h^{2}}{\cal O}_{h^{2}}{\cal O}_{h^{2}}\rangle_{0,1}+w({\cal O}_{h^{4}}{\cal O}_{1}|{\cal O}_{h^{2}})_{0,1}\no\\
&(=& 59021312+24850432\;\;)
\label{exam8}
\ea
This equality shares a fundamental structure with the two-point case. For $d=2$, Theorem \ref{gt} yields:
\ba
(1238948617930752&=)&\no\\
w({\cal O}_{h^{2}}{\cal O}_{h^{2}}|{\cal O}_{h^{2}})_{0,2}&=&\langle{\cal O}_{h^{2}}{\cal O}_{h^{2}}{\cal O}_{h^{2}}\rangle_{0,2}+\langle{\cal O}_{h^{2}}{\cal O}_{h^{2}}{\cal O}_{h^{2}}\rangle_{0,0}w({\cal O}_{h^{4}}{\cal O}_{1}|{\cal O}_{h^{2}})_{0,2}/k\no\\
&&+\langle{\cal O}_{h^{2}}{\cal O}_{h^{2}}{\cal O}_{h^{2}}{\cal O}_{h}\rangle_{0,1}w({\cal O}_{h^{5}}{\cal O}_{1})_{0,1}/k+\langle{\cal O}_{h^{2}}{\cal O}_{h^{2}}{\cal O}_{h^{2}}\rangle_{0,1}
w({\cal O}_{h^{4}}{\cal O}_{1}|{\cal O}_{h^{2}})_{0,1}/k\no\\
&(=&821654084851712+201251978293248\no\\
&&+59021312\cdot 4432896/8+59021312\cdot 24850432/8). 
\label{exam80} 
\ea 
In the right-hand side of the above equality, the first three terms exhibit characteristics similar to the two-point case, while the fourth term introduces a new effect specific to the multi-point scenario. We will discuss this example further in the next section.
\end{ex}

{\bf Acknowledgment} 
We would like to thank Prof. K. Hori for valuable discussions. 
Our research is partially supported by JSPS grants No. 22K03289 and No. 24H00182.  

\section{New Excess Intersections \\(\textcolor{black}{Erratum of page 463, lines 6--10 of the paper: Period integrals (Givental's $I$-function) of Calabi-Yau hypersurfaces in $CP^{N-1}$ as generating functions of intersection numbers on the moduli space of quasimaps \cite{JM}})  }
\label{intromp}
Let us begin by defining $Mp_{0,2|n}(N,d)$, which represents the highest dimensional stratum of $\widetilde{Mp}_{0,2|n}(N,d)$, the compactified moduli space of quasimaps from $CP^1$ with $2+n$ marked points to $CP^{N-1}$. For a detailed construction of $\widetilde{Mp}_{0,2|n}(N,d)$, please refer to our previous articles \cite{Jin1, JS}.

To describe a point in $Mp_{0,2|n}(N,d)$, we consider a vector-valued complex polynomial:
\ban
\sum_{j=0}^{d}{\bf a}_{j}s^{j}t^{d-j}\;\;\;\left({\bf a}_{j}\in {\bf C}^{N} (j=0,1,\cdots,d).\;{\bf a}_{0},{\bf a}_{d}\neq {\bf 0}\right).
\ean
Additionally, we have $z_j \in {\bf C}^\times$ for $j = 1, 2, \ldots, n$. Each $z_j \in {\bf C}^\times$ represents the point $(1:z_j) \in CP^1 - \{0, \infty\}$, where $0 = (1:0) \in CP^1$ and $\infty = (0:1) \in CP^1$. These data are represented by the tuple:
\ban
\left( ({\bf a}_0, {\bf a}_1, \ldots, {\bf a}_d), (z_1, z_2, \ldots, z_n) \right) \in ({\bf C}^N - \{{\bf 0}\}) \times ({\bf C}^N)^{d-1} \times ({\bf C}^N - \{{\bf 0}\}) \times ({\bf C}^\times)^n.
\ean
Next, we define a $({\bf C}^\times)^2$ action on $({\bf C}^N - \{{\bf 0}\}) \times ({\bf C}^N)^{d-1} \times ({\bf C}^N - \{{\bf 0}\}) \times ({\bf C}^\times)^n$ as follows:
\ba
&&G((\mu,\nu),\left(({\bf a}_{0},{\bf a}_{1},\cdots,{\bf a}_{d}),(z_{1},z_{2},\cdots,z_{n})\right))=\left((\mu{\bf a}_{0},\mu\nu{\bf a}_{1},\cdots,\mu\nu^{d}{\bf a}_{d}),(\nu z_{1},\nu z_{2},\cdots,\nu z_{n})\right)\no\\
&&\hspace{10.2cm}\left((\mu,\nu)\in ({\bf C}^{\times})^2\right). 
\label{twoC}
\ea
This action defines equivalence classes of elements in $({\bf C}^N - \{{\bf 0}\}) \times ({\bf C}^N)^{d-1} \times ({\bf C}^N - \{{\bf 0}\}) \times ({\bf C}^\times)^n$ based on the $({\bf C}^\times)^2$ orbits. We denote the equivalence class containing $\left( ({\bf a}_0, {\bf a}_1, \ldots, {\bf a}_d), (z_1, z_2, \ldots, z_n) \right)$ in the $({\bf C}^\times)^2$ orbit by $\left[ ({\bf a}_0, {\bf a}_1, \ldots, {\bf a}_d), (z_1, z_2, \ldots, z_n) \right]$. This equivalence class defines a point in $Mp_{0,2|n}(N,d)$. Note that $z_i$ and $z_j$ can coincide for $i \neq j$ in our definition.

In our previous work \cite{JS}, we defined the evaluation map $ev_m : \widetilde{Mp}_{0,2|n}(N,d) \rightarrow CP^{N-1}$ at the $m$-th marked point $(1:z_m) \in CP^1 - \{0, \infty\}$.
Let us take a point $\left[ ({\bf a}_0, {\bf a}_1, \ldots, {\bf a}_d), (z_1, z_2, \ldots, z_n) \right] \in Mp_{0,2|n}(N,d)$. Then, we consider the vector-valued polynomial $\sum_{j=0}^{d} {\bf a}_j s^j t^{d-j}$. It can be uniquely factored into the following form:
\ba
\sum_{j=0}^{d}{\bf a}_{j}s^{j}t^{d-j}=\left(\sum_{i=0}^{g}c_{i}s^{i}t^{g-i}\right)\cdot\left(\sum_{j=0}^{d-g}{\bf b}_{j}s^{j}t^{d-g-j}\right)\quad(c_{0},\cdots, c_{g}\in {\bf C},\;{\bf b}_{0},\cdots, {\bf b}_{d-g}\in {\bf C}^{N}), 
\label{factor0}
\ea  
where $\sum_{j=0}^{d-g} {\bf b}_j s^j t^{d-g-j}$ satisfies the condition:
\ba
\sum_{j=0}^{d-g}{\bf b}_{j}s^{j}t^{d-g-j}={\bf 0}\Longrightarrow (s,t)=(0,0).
\label{irr0}
\ea
With this setup, we ``may'' define $ev_m$ on $Mp_{0,2|n}(N,d)$ as follows:
\ba
ev_{m}([\left(({\bf a}_{0},{\bf a}_{1},\cdots,{\bf a}_{d}),(z_{1},z_{2},\cdots,z_{n})\right)]):=\pi_{N}\left(\sum_{j=0}^{d-g}{\bf b}_{j}(z_{m})^{d-g-j}\right),
\label{evm0}
\ea
where $\pi_N : ({\bf C}^N - \{{\bf 0}\}) \rightarrow CP^{N-1}$ is the projective equivalence map.

At first glance, the initial definition of the evaluation map $ev_m$ appears satisfactory, as it uniquely determines the image for any point in $Mp_{0,2|n}(N,d)$. However, it suffers from a significant discontinuity problem, as pointed out in \cite{hori}. 
\textcolor{black}{This discontinuity is the error in the definition of the evaluation map in \cite{JM}.}
Specifically, {\bf $ev_m$ can be discontinuous at points} satisfying:
\ba
[\left(({\bf a}_{0},{\bf a}_{1},\cdots,{\bf a}_{d}),(z_{1},z_{2},\cdots,z_{n})\right)]\in Mp_{0,2|n}(N,d) \;\;\mbox{which satisfies}\;\; \sum_{j=0}^{d}{\bf a}_{j}(z_{m})^{d-j}={\bf 0}.
\ea 
This discontinuity arises when the non-negative integer $g$ in the factorization (\ref{factor0}) is a positive integer and 
\ba
\sum_{i=0}^{g}c_{i}(z_{m})^{g-i}=0.
\label{c0}
\ea
 To illustrate this, let us consider a vector-valued polynomial $\sum_{j=0}^{d}{\bf e}_{j}s^{j}t^{d-j}$ that satisfies the condition:
\ba
\sum_{j=0}^{d}{\bf e}_{j}s^{j}t^{d-j}={\bf 0}\;\Longrightarrow\; (s,t)=(0,0). 
\label{disc}
\ea 
Now, consider the following deformation of the original polynomial:
\ba
\varphi_{\varepsilon}(s,t)&:=&\sum_{j=0}^{d}{\bf a}_{j}s^{j}t^{d-j}+\varepsilon\left(\sum_{j=0}^{d}{\bf e}_{j}s^{j}t^{d-j}\right)\no\\
                                 &=&\left(\sum_{i=0}^{g}c_{i}s^{i}t^{g-i}\right)\cdot\left(\sum_{j=0}^{d-g}{\bf b}_{j}s^{j}t^{d-g-j}\right)+\varepsilon\left(\sum_{j=0}^{d}{\bf e}_{j}s^{j}t^{d-j}\right).
\label{defor}
\ea   
When $\varepsilon=0$, $ev_m([\left(({\bf a}_{0},{\bf a}_{1},\cdots,{\bf a}_{d}),(z_{1},z_{2},\cdots,z_{n})\right)])$ is given by $[\sum_{j=0}^{d-g}{\bf b}_{j}(z_{m})^{d-g-j}]$ as defined in (\ref{evm0}). However, for any non-zero $\varepsilon$, the image abruptly changes to $\left[\varepsilon\left(\sum_{j=0}^{d}{\bf e}_{j}(z_{m})^{d-j}\right)\right]=[\sum_{j=0}^{d}{\bf e}_{j}(z_{m})^{d-j}]$ because (\ref{c0}) holds. Consequently, the image jumps from $[\sum_{j=0}^{d-g}{\bf b}_{j}(z_{m})^{d-g-j}]$ to $[\sum_{j=0}^{d}{\bf e}_{j}(z_{m})^{d-j}]$, indicating that $ev_m$ is discontinuous at points satisfying the factorization condition (\ref{disc}) with positive $g$.  

To resolve this discontinuity issue, we adopt the approach proposed by Ciocan-Fontanine and Kim \cite{CK2}, which involves:
\begin{itemize}
\item[(i)]{Changing the target domain of $ev_m$ from $CP^{N-1}$ to the algebraic stack ${\bf C}^{N}/{\bf C}^{\times}$.}
\end{itemize}
This stack may be viewed as $CP^{N-1}\cup\{{\bf 0}\}$ as a quotient set. However, with the quotient topology, it is not a Hausdorff space, due to the following property:
\begin{itemize}
\item[(ii)]{Every point $p\in CP^{N-1}$ is infinitesimally close to ${\bf 0}$, meaning any open set containing ${\bf 0}$ also contains the entire $CP^{N-1}$.}
\end{itemize}
We therefore modify the definition of $ev_m$ as follows:
\begin{defi}
\ban
ev_{m}([\left(({\bf a}_{0},{\bf a}_{1},\cdots,{\bf a}_{d}),(z_{1},z_{2},\cdots,z_{n})\right)]):=\tilde{\pi}_{N}\left(\sum_{j=0}^{d}{\bf a}_{j}(z_{m})^{d-j}\right),
\label{evm}
\ean
where $\tilde{\pi}_{N}:{\bf C}^{N}\rightarrow {\bf C}^{N}/{\bf C}^{\times}$ is the quotient map.
\label{newev}
\end{defi}
Now, let us reconsider the image of the point corresponding to $\varphi_{\varepsilon}(s,t)$ in (\ref{defor}) with $\sum_{i=0}^{g}c_{i}(z_{m})^{g-i}=0$. When $\varepsilon=0$, the image is ${\bf 0}$. When $\varepsilon\neq 0$, it becomes $[\sum_{j=0}^{d}{\bf e}_{j}(z_{m})^{d-j}]$. Due to property (ii), ${\bf 0}$ is infinitesimally close to $[\sum_{j=0}^{d}{\bf e}_{j}(z_{m})^{d-j}]$, ensuring that $ev_m$ is continuous at this point.

\textcolor{black}{In sum, we write down the Erratum of \cite{JM} in order to resolve the discontinuity problem of the evaluation map.}
\begin{err}{\bf (\cite{JM})}
The definition of $ev_{1}:Mp_{0,2|1}(N,d)\rightarrow CP^{N-1}$ on page 463, 
lines 6–10 of \cite{JM} (stated as (\ref{evm0}) in the present paper) should be replaced with Definition \ref{newev} for $m=1$.
\end{err}

With these preparations, we now discuss the novel excess intersections that become relevant in the proof of Theorem \ref{gt}. As demonstrated in \cite{Jin6}, the generalized mirror transformation for virtual structure constants with two operator insertions arises from excess intersections. In essence, we can interpret the mirror map as {\bf the generating function of these excess intersections}. These intersections primarily manifest because the dimension of the moduli space either decreases or remains invariant as the degree $d$ increases, particularly when the degree of the hypersurface $k$ is no less than the positive integer $N$ of $CP^{N-1}$. However, Theorem \ref{gt} reveals that the generalized mirror transformation is non-trivial even when $N$ exceeds $k$. Consequently, we must identify new excess intersections that account for the non-trivial terms in the theorem. These intersections can indeed be found by considering the revised definition of $ev_m$ and the characteristic property (ii) of ${\bf C}^{N}/{\bf C}^{\times}$.

Let us fix $g$ to a positive integer no greater than $d$ and examine vector-valued polynomials factorized as: 
\ba
\varphi(s,t)&=&\sum_{j=0}^{d}{\bf a}_{j}s^{j}t^{d-j}\no\\
              &=&\left(\sum_{i=0}^{g}c_{i}s^{i}t^{g-i}\right)\cdot\left(\sum_{j=0}^{d-g}{\bf b}_{j}s^{j}t^{d-g-j}\right),
\label{fac2}              
\ea
where $\sum_{j=0}^{d-g}{\bf b}_{j}s^{j}t^{d-g}$ satisfies the condition (\ref{irr0}). We further factorize $\sum_{i=0}^{g}c_{i}s^{i}t^{g-i}$ into the form:
\ba
\sum_{i=0}^{g}c_{i}s^{i}t^{g-i}=c_{0}\prod_{i=1}^{l}(t-w_{i}s)^{g_{i}}\;\;\quad \left((g_{1},g_{2},\cdots,g_{l})=\sigma_{g}\in P_{g}^{l}\right),
\label{cfac}
\ea
where the $w_i$'s are $l$ distinct elements of ${\bf C}^{\times}$. With this setup, we consider the locus ${\cal L}_{(m,i)} \subset \widetilde{Mp}_{0,2|n}(N,d)$ where the $m$-th marked point $z_m$ coincides with $w_i$. We then analyze the contribution of the locus under consideration to the multi-point virtual structure constant:
\ba
w({\cal O}_{h^{a}}{\cal O}_{h^{b}}|\prod_{j=1}^{n}{\cal O}_{h^{c_{j}}})_{0,d}=\int_{\widetilde{Mp}_{0,2|n}(N,d)}
c_{top}(\widetilde{\cal E}_{d}^{k})\wedge ev_{0}^{*}(h^a)\wedge ev_{\infty}^{*}(h^b)\wedge\left(\mathop{\bigwedge}_{j=1}^{n}ev_{j}^{*}(h^{c_{j}})\right).
\label{w2}
\ea  
Given the factorization condition (\ref{fac2}), the number of degrees of freedom arising from $\sum_{j=0}^{d-g}{\bf b}_{j}s^{j}t^{d-g-j}$ and $\sum_{i=0}^{g}c_{i}s^{i}t^{g-i}$ is $N-3+(N-k)(d-g)+l$\footnote{We must account for the constraint that the image of the polynomial map lies within the hypersurface $M_{N}^{k}$, reflected by the insertion of $c_{top}(\widetilde{\cal E}_{d}^{k})$ in (\ref{defw}). This constraint decreases the degrees of freedom of $\sum_{j=0}^{d-g}{\bf b}_{j}s^{j}t^{d-g-j}$, $N(d-g+1)-2=N-2+N(d-g)$ ($-2$ comes from the $({\bf C}^{\times})^2$ action in (\ref{twoC})) into $N-3+(N-k)(d-g)$. }. Incorporating the degrees of freedom from the marked points other than $z_m$, the total dimension of the locus ${\cal L}_{(m,i)}$ is:
\ban
N-3+(N-k)(d-g)+l+n-1.
\ean 
Applying the condition (\ref{wsel}), the number of degrees of freedom of ${\cal L}_{(m,i)}$ minus the constraints imposed by operator insertions is:
\ba
&&N-3+(N-k)(d-g)+l+n-1-\left(a+b+\sum_{j=1}^{n}c_{j}\right)\no\\
&=&N-3+(N-k)(d-g)+l+n-1-\left(N-3+(N-k)d+n\right)\no\\
&=&l-1-(N-k)g\no\\
&=&-(g-l)-1-(N-k-1)g.
\label{cntdim}
\ea
Since $g\geq l$, the integer in the last line is consistently negative if $N-k-1\geq 0$. Assuming $ev_{m}^{*}(h^{c_{m}})$ is an ordinary cohomology element of (complex) degree $c_m$, we conclude that the locus under consideration does not contribute to the intersection number (\ref{w2}) when $N-k-1\geq 0$.

At this juncture, we must emphasize the condition where the marked point $z_m$ coincides with $w_i$, expressed as:
\ba
\tilde{\pi}_{N}\left(\varphi(1,z_{m})\right)={\bf 0} \in {\bf C}^{N}/{\bf C}^{\times}\Longleftrightarrow ev_{m}([\left(({\bf a}_{0},{\bf a}_{1},\cdots,{\bf a}_{d}),(z_{1},z_{2},\cdots,z_{n})\right)])
={\bf 0} \in {\bf C}^{N}/{\bf C}^{\times}.
\label{evme}
\ea
The right-hand side of (\ref{w2}) indicates that $w({\cal O}_{h^{a}}{\cal O}_{h^{b}}|\prod_{j=1}^{n}{\cal O}_{h^{c_{j}}})_{0,d}$ generically counts the number of quasimaps within the following subset of $\widetilde{Mp}_{0,2|n}(N,d)$:
\ban
PD(c_{top}({\cal E}_{d}^{k}))\cap ev_{0}^{-1}(PD(h^{a}))\cap ev_{\infty}^{-1}(PD(h^{b}))\cap\left(\mathop{\bigcap}_{j=1}^{m}ev_{j}^{-1}(PD(h^{c_{j}}))\right),
\ean
where $PD(*)$ denotes the subvariety Poincar\'{e} dual to the cohomology element $*$. Let us examine $ev_{m}^{-1}(PD(h^{c_{m}}))$ under the conditions of (\ref{evme}). Recall that  $PD(h^{c_{m}})$ is a linear
space of codim $c_{m}$  in $CP^{N-1}$, and the fact that  ${\bf 0}\in {\bf C}^{N}/{\bf C}^{\times}$ is infinitesimally close to (or indistinguishable in the quotient topology from) any point in $CP^{N-1}$.
The condition (\ref{evme}) implies
\ba
ev_{m}^{-1}(PD(h^{c_{m}}))=ev_{m}^{-1}({\bf C}^{N}/{\bf C}^{\times})\quad\mbox{on ${\cal L}_{(m,i)}$},
\ea
meaning {\bf $ev_{m}^{*}(h^{c_{m}})$ imposes no constraints on ${\cal L}_{(m,i)}$}.
This conclusion is also derived from elementary consideration. Since $PD(h^{c_{m}})$ is a linear space, restriction to $ev_{m}^{-1}(PD(h^{c_{m}}))$ is given by linear equations of homogeneous coordinates 
of the image $ev_{m}([\left(({\bf a}_{0},{\bf a}_{1},\cdots,{\bf a}_{d}),(z_{1},z_{2},\cdots,z_{n})\right)])$.
But the condition (\ref{evme}) turns these linear equations into $0=0$, i.e., trivial equations. 
Consequently, the dimensional count of (\ref{cntdim}) must be adjusted as follows:
\ba
&&N-3+(N-k)(d-g)+l+n-1-\left(a+b+\sum_{j\neq m}^{}c_{j}\right)\no\\
&=&N-3+(N-k)(d-g)+l+n-1-\left(N-3+(N-k)d+n-c_{m}\right)\no\\
&=&c_{m}+l-1-(N-k)g\no\\
&=&(c_{m}-1)-(g-l)-(N-k-1)g.
\label{cntdim2}
\ea
This adjustment allows the last line to yield a non-negative integer. This is a case where the locus ${\cal L}_{(m,i)}$ contributes to the intersection number (\ref{w2}) as new excess intersections even when $N-k-1\geq 0$. More generally, we can consider scenarios where $p$ marked points $z_{m_{1}},z_{m_{2}},\cdots,z_{m_{p}}\;(1\leq m_{1}<m_{2}<\cdots<m_{p}\leq n)$ coincide with elements of $\{w_{1},w_{2},\cdots,w_{l}\}$. In such cases, the dimensional count for the new excess intersection is modified to:
\ba
&&N-3+(N-k)(d-g)+l+n-p-\left(a+b+\sum_{j\notin \{m_{1},\cdots,m_{p}\}}^{}c_{j}\right)\no\\
&=&\sum_{j=1}^{p}(c_{m_{j}}-1)-(g-l)-(N-k-1)g.
\label{cntdim3}
\ea
\begin{ex}
The generalized mirror transformation for $w({\cal O}_{h}{\cal O}_{h^{2}}|{\cal O}_{h^{2}})_{0,1}$ in (\ref{exam2}) is trivial, as demonstrated by setting $c_{1}=2,g=l=1,N=4,k=1$ in (\ref{cntdim2}):
\ban
2-1-(1-1)-(4-1-1)\cdot 1=-1<0.
\ean
However, for $w({\cal O}_{h}{\cal O}_{h}|{\cal O}_{h^{2}}{\cal O}_{h^{2}})_{0,1}$ in (\ref{exam20}), we obtain from (\ref{cntdim3}) with $c_{1}=c_{2}=2, g=l=1,N=4,k=1$:
\ban
(2-1)+(2-1)-(1-1)-(4-1-1)\cdot 1=0.
\ean
Thus, the term $w({\cal O}_{h^2}{\cal O}_{1}|{\cal O}_{h^{2}}{\cal O}_{h^{2}})_{0,1}$ arises from the new excess intersection. 
\end{ex}
\begin{ex}
For $w({\cal O}_{h^2}{\cal O}_{h^2}|{\cal O}_{h^2})_{0,1}$ in (\ref{exam8}), the term $w({\cal O}_{h^4}{\cal O}_{1}|{\cal O}_{h^2})_{0,1}$ emerges from (\ref{cntdim2}) with $c_{1}=2, g=l=1,N=k=8$:
\ban
(2-1)-(1-1)-(8-8-1)\cdot 1=2\geq 0.
\ean 
In the case of $w({\cal O}_{h^2}{\cal O}_{h^{2}}|{\cal O}_{h^{2}})_{0,2}$ in (\ref{exam80}), the term $\langle{\cal O}_{h^{2}}{\cal O}_{h^2}{\cal O}_{h^2}\rangle_{0,0}w({\cal O}_{h^{4}}{\cal O}_{1}|{\cal O}_{h^2})_{0,2}/k$ appears from (\ref{cntdim2}) with $c_1=2,g=2,l=1,N=k=8$:
\ban
2-1-(2-1)-(8-8-1)\cdot 2=2\geq 0.
\ean
$\langle{\cal O}_{h^{2}}{\cal O}_{h^{2}}{\cal O}_{h^{2}}{\cal O}_{h}\rangle_{0,1}w({\cal O}_{h^5}{\cal O}_{1})_{0,1}/k$ results from the excess intersection already discussed in \cite{Jin6}, with the following count:
\ban
-(g-l)-(N-k-1)g=-(1-1)-(8-8-1)\cdot 1=1\geq 0.
\ean
The term $\langle{\cal O}_{h^{2}}{\cal O}_{h^{2}}{\cal O}_{h^{2}}\rangle_{0,1}w({\cal O}_{h^{4}}{\cal O}_{1}|{\cal O}_{h^{2}})_{0,1}/k$ arises from the new excess intersection (\ref{cntdim2}) with $c_1=2,g=l=1,N=k=8$:
\ban
2-1-(1-1)-(8-8-1)\cdot 1 =2\geq 0.
\ean
\end{ex}

\section{Proof of Theorem \ref{gt}}
The term $w({\cal O}_{h^{a}}{\cal O}_{h^{b}}
|\prod_{j=2}^{N-2}({\cal O}_{h^{j}})^{n_{j}})_{0,d}$ on the left-hand side of (\ref{main}) is defined as an intersection number on the moduli space $\widetilde{Mp}_{0,2|\sum_{j=2}^{N-2}n_{j}}(N,d)$.
In this proof, we focus on the highest dimensional stratum $Mp_{0,2|\sum_{j=2}^{N-2}n_{j}}(N,d).$\footnote{See the last part of the proof of the two-point version of Theorem \ref{gt} \cite{Jin6}. 
We will also explain the reason in the last part of this paper.}
We introduce the notational setting for $Mp_{0,2|\sum_{j=2}^{N-2}n_{j}}(N,d)$: 
\ba
&&Mp_{0,2|\sum_{j=2}^{N-2}n_{j}}(N,d)=\no\\
&&\{\left(({\bf a}_{0},{\bf a}_{1},\cdots,{\bf a}_{d}),(z_{(2,1)},\cdots,z_{(2,n_{2})},\cdots,z_{(N-2,1)},\cdots,z_{(N-2,n_{N-2})})\right) \;|\;{\bf a}_{i}\in{\bf C}^{N},{\bf a}_{0},{\bf a}_{d}\neq {\bf 0},z_{(j,i)}\in {\bf C}^{\times}\;\}
/({\bf C}^{\times})^{2}, \no\\
\ea
where the two ${\bf C}^{\times}$ actions are given by:
\ba
&&\lambda\cdot\left(({\bf a}_{0},{\bf a}_{1},\cdots,{\bf a}_{d}),(z_{(2,1)},\cdots,z_{(2,n_{2})},\cdots,z_{(N-2,1)},\cdots,z_{(N-2,n_{N-2})})\right)=\no\\ 
&&\left((\lambda {\bf a}_{0},\lambda{\bf a}_{1},\cdots,\lambda{\bf a}_{d}),(z_{(2,1)},\cdots,z_{(2,n_{2})},\cdots,z_{(N-2,1)},\cdots,z_{(N-2,n_{N-2})})\right) ,\no\\
&&\nu\cdot \left(({\bf a}_{0},{\bf a}_{1},\cdots,{\bf a}_{d}),(z_{(2,1)},\cdots,z_{(2,n_{2})},\cdots,z_{(N-2,1)},\cdots,z_{(N-2,n_{N-2})})\right)=\no\\ 
&&\left(({\bf a}_{0},\nu{\bf a}_{1},\nu^{2}{\bf a}_{2},\cdots,\nu^{d}{\bf a}_{d}),(\nu z_{(2,1)},\cdots,\nu z_{(2,n_{2})},\cdots,\nu z_{(N-2,1)},\cdots,\nu z_{(N-2,n_{N-2})})\right),\no\\
&&\hspace{10.8cm}(\lambda,\nu \in {\bf C}^{\times}).
\label{equi}
\ea
The points $(1:z_{(j,i)})\in (CP^{1}-\{0,\infty\})$ for $(i=1,\cdots, n_{j})$ are the marked points where the operator ${\cal O}_{h^{j}}$ is inserted.
Note that the $z_{(j,i)}$'s can move freely within ${\bf C}^{\times}$, and they may coincide. 
The tuple $({\bf a}_{0},{\bf a}_{1},\cdots,{\bf a}_{d-1},{\bf a}_{d})$ represents a polynomial map $\varphi(s:t)=\tilde{\pi}_{N}\left(\sum_{j=0}^{d}{\bf a}_{j}s^{j}t^{d-j}\right)$ from $CP^{1}$ to 
${\bf C}^{N}/{\bf C}^{\times}$, modulo the $\nu$-action, which is introduced from the ${\bf C}^{\times}$ action on $CP^{1}$ that fixes $0=(0:1), \infty=(1:0)\in CP^{1}$:
\ba
(s:t)\rightarrow (\nu s:t).
\label{p1}
\ea
As was discussed in the previous section, $\sum_{j=0}^{d}{\bf a}_{j}s^{j}t^{d-j}$ can be uniquely factorized into the following form:
\ba
\sum_{j=0}^{d}{\bf a}_{j}s^{j}t^{d-j}=c_{0}\biggl(\prod_{j=1}^{l}(t-w_{j}s)^{g_{j}}\biggr)
\cdot\biggl(\sum_{j=0}^{d-g}{\bf b}_{j}s^{j}t^{d-g-j}\biggr),
\label{factor}
\ea
where the following condition holds:
\ban
\sum_{j=0}^{d-g}{\bf b}_{j}s^{j}t^{d-g-j}={\bf 0}\Longrightarrow (s,t)=(0,0),
\label{irr}
\ean
or, equivalently, sending $(s,t)$ to $\pi_{N}(\sum_{j=0}^{d-g}{\bf b}_{j}s^{j}t^{d-g-j})$ defines a well-defined map from $CP^{1}$ to $CP^{N-1}$.
 
Since ${\bf a}_{0},{\bf a}_{d}\neq {\bf 0}$, it follows that $w_{j}$ belongs to ${\bf C}^{\times}$, which means $(1:w_{j})$ never coincides with $0$ or $\infty$ in $CP^{1}$. These distinct points $(1:w_{j})\;(j=1,\cdots,l(\sigma_{g}))$ also represent points where $\varphi(s:t)$ evaluates to ${\bf 0}$. Given that $[({\bf a}_{0},{\bf a}_{1},\cdots,{\bf a}_{d-1},{\bf a}_{d})]$ represents $\varphi(s:t)$ modulo the ${\bf C}^{\times}$ action defined in (\ref{p1}), the configuration of the $(1:w_{j})$'s must be considered modulo this ${\bf C}^{\times}$ action. For brevity, we will henceforth not distinguish between $w_{j}\in {\bf C}^{\times}$ and $(1:w_{j})\in (CP^{1}-\{0,\infty\})$, nor between $z_{(j,i)}\in {\bf C}^{\times}$ and $(1:z_{(j,i)})\in (CP^{1}-\{0,\infty\})$.

It is evident that the factorization in (\ref{factor}) remains invariant under permutations of $w_{j}$'s with subscripts $j$ that share the same value of $g_{j}$. Based on these observations, we obtain the following decomposition of $Mp_{0,2|\sum_{j=2}^{N-2}n_{j}}(N,d)$:
\ba
Mp_{0,2|\sum_{j=2}^{N-2}n_{j}}(N,d)=\coprod_{g=0}^{d}\coprod_{l=1}^{g}\coprod_{\sigma_{g}\in P_{g}^{l}}
M_{0,2+l|\sum_{j=2}^{N-2}n_{j}}(CP^{N-1},d-g,\sigma_{g}).
\ea
This decomposition includes $M_{0,2+l|\sum_{j=2}^{N-2}n_{j}}(CP^{N-1},d-g,\sigma_{g})$, which denotes the uncompactified moduli space of degree $d-g$ holomorphic maps from $CP^{1}$ to $CP^{N-1}$. This space features $2+l$ distinct marked points and $\sum_{j=2}^{N-2}n_{j}$ marked points free to move within $CP^{1}-\{0,\infty\}$, modulo the action of $\prod_{i=1}^{g}\mbox{Sym}(\mbox{mul}(i,\sigma_{g}))$ on the set $\{w_{1}, w_{2},\cdots, w_{l}\}$.
\ba
&&M_{0,2+l|\sum_{j=2}^{N-2}n_{j}}(CP^{N-1},d-g,\sigma_{g})\no\\
&&:=\{[\bigl(\tilde{\pi}_{N}(\sum_{j=0}^{d-g}{\bf b}_{j}s^{j}t^{d-g-j}),(0,\infty,w_{1},w_{2},\cdots,w_{l}), (z_{(2,1)},\cdots,z_{(N-2,n_{N-2})})\bigr)]\}/\biggl(\prod_{i=1}^{g}\mbox{Sym}(\mbox{mul}(i,\sigma_{g}))\biggr).\no\\
\ea 
In this definition, the tuple $\bigl(\pi_{N}(\sum_{j=0}^{d-g}{\bf b}_{j}s^{j}t^{d-g-j}),(0,\infty,w_{1},w_{2},\cdots,w_{l}),(z_{(2,1)},\cdots,z_{(N-2,n_{N-2})})\bigr)$ is considered modulo the ${\bf C}^{\times}$ action 
(\ref{p1}) on $CP^{1}$ , and the symmetric group $\mbox{Sym}(\mbox{mul}(i,\sigma_{g}))$ permutes the $w_{j}$'s for which $g_{j}=i$.

As observed in the discussion given in first half of Section 3 of \cite{Jin6}, ill-definedness of the image of the polynomial map $\varphi$ at the points $w_{j}$ leads to the excess intersections in the two-pointed case. These remain relevant in our current analysis. However, we must also account for the new excess intersections introduced in the preceding section. To this end, we define a decomposition of the subscript set for the $z_{(*,*)}$'s as follows:
\ba
&&\{(j,1),(j,2),\cdots,(j,n_{j})\}=\mathop{\coprod}_{p=0}^{l}\{(j,h_{1}^{p}),(j,h_{2}^{p}),\cdots,(j,h_{m_{p}^{j}}^{p})\}
=:\mathop{\coprod}_{p=0}^{l}W_{(j,p)}\no\\ 
&& (j=2,3,\cdots, N-2,\;\;m_{p}^{j}\geq 0,\;\;\;m_{0}^{j}+\sum_{p=1}^{l}m_{p}^{j}=n_{j}\;\;\;\;1\leq h_{1}^{p}<h_{2}^{p}<\cdots<h_{m_{p}^{j}}^{p}\leq n_{j},\;\;).
\ea
Under this decomposition, we impose the following restrictions on the positions of the $z_{(j,i)}$'s:
\ba
&&z_{(j,h^{0}_{i})}\notin \{w_{1},w_{2},\cdots, w_{l}\}\;\;(i=1,2,\cdots, m_{0}^{j}),\no\\
&&z_{(j,h^{p}_{i})}=w_{p}\;\;(p=1,2,\cdots,l,\;\;i=1,2,\cdots, m_{p}^{j}).
\label{rest}
\ea
For fixed positive integers $m_{a}^{j}$'s $(a=0,1,\cdots,l, j=2,3,\cdots,N-2)$, the number of distinct decompositions is given by: 
\ba
\prod_{j=2}^{N-2}\frac{n_{j}!}{m_{0}^{j}!\left(\prod_{p=1}^{l}m_{p}^{j}!\right)}.
\label{comb}
\ea 
Let us fix a decomposition type $\mathop{\coprod}_{j=2}^{N-2}\mathop{\coprod}_{p=0}^{l}W_{(j,p)}$ and denote by $M_{0,2+l| \mathop{\mathop{\coprod}_{j=2}^{N-2}\coprod}_{p=0}^{l}W_{(j,p)}}(CP^{N-1},d-g,\sigma_{g})$ the locus within $M_{0,2+l|\sum_{j=2}^{N-2}n_{j}}(CP^{N-1},d-g,\sigma_{g})$ where the $z_{(*,*)}$'s satisfy the restrictions of (\ref{rest}). 
To evaluate the contribution from the excess intersection arising from this locus, we introduce a 
\textcolor{black}{``perturbation space''} defined as follows \footnote{The logic behind the introduction of this perturbation space mirrors that used in \cite{Jin6}. See Appendix A of \cite{Jin6}.}:
\ba
&&M^{pert.}_{0,2+l| \mathop{\coprod}_{j=2}^{N-2}\mathop{\coprod}_{p=0}^{l}W_{(j,p)}}(CP^{N-1},d-g,\sigma_{g} )\no\\
&&:=\left(M_{0,2+l| \mathop{\coprod}_{j=2}^{N-2}W_{(j,0)}}(CP^{N-1},d-g)\biggl(\prod_{p=1}^{l}(\mathop{\times}_{CP^{N-1}}
\widetilde{Mp}_{0,2|\sum_{j=2}^{N-2}m_{p}^{j}}(N,g_{p}))\biggr)\right)/\biggl(\prod_{i=1}^{g}\mbox{Sym}(\mbox{mul}(i,\sigma_{g}))\biggr).\no\\
\label{pertsp}
\ea 
Here, $M_{0,2+l \mathop{\coprod}_{j=2}^{N-2}W_{(j,0)}}(CP^{N-1},d-g)$ represents the uncompactified moduli space of holomorphic maps from $CP^{1}$ with $2+l$ distinct marked points 
and the marked points $z_{(*,*)}$'s whose subscripts belong to $\mathop{\coprod}_{j=2}^{N-2}W_{(j,0)}$.
The $p$-th fiber product in (\ref{pertsp}) is defined through the following diagram:
\ba
&&M_{0,2+l(\sigma_{g})| \mathop{\coprod}_{j=2}^{N-2}W_{(j,0)}}(CP^{N-1},d-g)\stackrel{ev_{w_{p}}}{\longrightarrow}CP^{N-1}\stackrel{\tiny{ev_{0}}}{\longleftarrow}
\widetilde{Mp}_{0,2|\sum_{j=2}^{N-2}m_{p}^{j}}(N,g_{p}))\;\;\;\; ( p=1,2,\cdots,l).\no\\
\label{dgm1}
\ea
In the above diagram, $ev_{w_{p}}$ is the evaluation map of $\tilde{\pi}_{N}(\sum_{j=0}^{d-g}{\bf b}_{j}s^{j}t^{d-g-j})$ at $w_{p}$.
The symmetric group $\prod_{i=1}^{g}\mbox{Sym}(\mbox{mul}(i,\sigma_{g}))$ permutes the $w_{i}$'s along with $\mathop{\times}_{CP^{N-1}}\widetilde{Mp}_{0,2| \sum_{j=2}^{N-2}m_{i}^{j}}(N,g_{i})$ in a manner consistent with the definition of $M_{0,2+l|\sum_{j=2}^{N-2}n_{j}}(CP^{N-1},d-g,\sigma_{g})$.

The corresponding contribution is then given by:
\ba
&&\int_{M^{pert.}_{0,2+l(\sigma_{g})|\mathop{\coprod}_{j=2}^{N-2} \mathop{\coprod}_{p=0}^{l}W_{(j,p)}}(CP^{N-1},d-g,\sigma_{g} ))}c_{top}(\widetilde{\cal E}^{k}_{d})\wedge ev_{0}^{*}(h^{a})
\wedge ev_{\infty}^{*}(h^b)\wedge\left(\mathop{\bigwedge}_{\substack{j=2,\cdots,N-2\\i=1,2,\cdots,n_{j}}}ev_{(j,i)}^{*}(h^{j})\right).\no\\
&=&
\sum_{e_{1},\cdots,e_{l}=0}^{N-2}
\langle{\cal O}_{h^{a}}{\cal O}_{h^{b}}\prod_{j=2}^{N-2}({\cal O}_{h^{j}})^{m_{0}^{j}}\prod_{p=1}^{l}{\cal O}_{h^{e_{p}}}\rangle_{0,d-g}
\biggl(\prod_{p=1}^{l}\frac{w({\cal O}_{h^{N-2-e_{p}}}{\cal O}_{1}|\prod_{j=2}^{l}({\cal O}_{h^{j}})^{m_{p}^{j}})_{0,g_{p}}}{k}\biggr)
\biggl(\prod_{i=1}^{g}\frac{1}{(\mbox{mul}(i,\sigma_{g}))!}\biggr), \no\\
\label{cont}  
\ea
This follows because $w({\cal O}_{h^{N-2-e_{i}}}{\cal O}_{1}|\prod_{j=2}^{l}({\cal O}_{h^{j}})^{m_{p}^{j}})_{0,g_{p}}$ (which includes the insertion of ``1'', or no insertion) is non-vanishing. 
In deriving this formula, we have used the splitting axiom \cite{KM} (2.2.6. page 529)   for the fiber product, i.e., the fiber 
product $\mathop{\times}_{CP^{N-1}}$ under the diagram:
\ban
 M_{0,2+l|\mathop{\coprod}_{j=2}^{N-2}W_{(j,0)}}(CP^{N-1},d-g)\stackrel{ev_{p}}{\longrightarrow}
CP^{N-1}\stackrel{ev_{0}}{\longleftarrow} \widetilde{Mp}_{0,2}(N,g_{p})
\ean
corresponds to,
\ban 
\frac{1}{k}\sum_{e_{p}=0}^{N-2}ev_{w_{p}}^{*}(h^{e_{p}})\wedge ev_{0}^{*}(h^{N-2-e_{p}}).
\ean
This correspondence arises from the description of the Poincar\'{e} dual of the diagonal $\Delta$ in $CP^{N-1}\times  CP^{N-1}$ (specifically, $M_{N}^{k}\times M_{N}^{k}$) using cohomology elements.

The term $\langle{\cal O}_{h^{a}}{\cal O}_{h^{b}}\prod_{j=2}^{N-2}({\cal O}_{h^{j}})^{m_{0}^{j}}\prod_{i=1}^{l}{\cal O}_{h^{e_{i}}}\rangle_{0,d-g}$ in (\ref{cont}) represents the Gromov-Witten invariant of the hypersurface, obtained by integration over $M_{0,2+l| \mathop{\coprod}_{j=2}^{N-2}W_{(j,0)}}(CP^{N-1},d-g)$.
While the moduli space $M_{0,2+l| \mathop{\coprod}_{j=2}^{N-2}W_{(j,0)}}(CP^{N-1},d-g)$ is uncompactified and the marked points $z_{(j,i)}$'s can coincide, the boundary components added during compactification by stable maps, and the coincidence loci of the $z_{(j,i)}$'s, do not contribute to the intersection numbers due to their positive complex codimensions. 
The factor $\displaystyle{\prod_{i=1}^{g}\frac{1}{(\mbox{mul}(i,\sigma_{g}))!}}$  on the right-hand side of (\ref{cont}) accounts for the division by the group $\prod_{i=1}^{g}\mbox{Sym}(\mbox{mul}(i,\sigma_{g}))$.

In the $g=0$ case, the right-hand side of (\ref{cont}) simplifies to $\langle{\cal O}_{h^{a}}{\cal O}_{h^{b}}\prod_{j=2}^{N-2}({\cal O}_{h^{j}})^{n_{j}}\rangle_{0,d}$.
In the $g=d$ case, $\langle{\cal O}_{h^{a}}{\cal O}_{h^{b}}\prod_{j=2}^{N-2}({\cal O}_{h^{j}})^{m_{0}^{j}}\prod_{i=1}^{l}{\cal O}_{h^{e_{i}}}\rangle_{0,0}$ vanishes if $l>1$, or if there exists a $j\in\{2,3,\cdots,N-2\}$ such that $m_{0}^{j}>0$ \cite{KM} (2.2.5. page 529). 
If $l=1$ and $m_{0}^{2}=m_{0}^{3}=\cdots=m_{0}^{N-2}=0$, we have:
\ba
 \langle{\cal O}_{h^a}{\cal O}_{h^{b}}{\cal O}_{h^{e_{1}}}\rangle_{0,0}=k\cdot\delta_{e_{1},N-2-a-b}.
\ea
Consequently, the right-hand side of (\ref{cont}) becomes $w({\cal O}_{h^{a+b}}{\cal O}_{1}|\prod_{j=2}^{N-2}({\cal O}_{h^{j}})^{n_{j}})_{0,d}$.
We can also introduce perturbation spaces for the lower-dimensional strata of $\widetilde{Mp}_{0,2|\sum_{j=2}^{N-2}n_{j}}(N,d)$, but these are irrelevant to the intersection number due to their positive codimensions.

By multiplying the contribution (\ref{cont}) by the combinatorial factor given in (\ref{comb}), and summing over all possible choices of $g$'s, $l$'s, $\sigma_{g}$'s and $m_{p}^{j}$'s, we obtain the right-hand side of Theorem \ref{gt}.        \hspace{2.2cm} $\Box$

 \vspace{1em}
 \noindent
 {\bf Data Availability Statement}  This article has no associated data. \\

\end{document}